\documentclass[11pt,reqno]{amsart}
\usepackage{amsmath,amsfonts,amssymb,amsthm}
\usepackage{psfig}
\newlength{\figboxwidth}
\newcommand{\bdry}{\partial}

\def\@ifundefined#1#2#3%
  {\expandafter\ifx\csname#1\endcsname\relax#2\else#3\fi}

\@ifundefined{theoremstyle}{
}{
\theoremstyle{plain} 
}
\newtheorem{theorem}{Theorem}[section]

\newtheorem{proposition}[theorem]{Proposition}
\newtheorem{lemma}[theorem]{Lemma}

\@ifundefined{theoremstyle}{
}{
\theoremstyle{definition} 
}

\newlength{\halfbls}

\newcommand\NN{{\mathcal N}}

\catcode`\@=11
\long\def\@savemarbox#1#2{\global\setbox#1\vtop{\hsize\marginparwidth 
  \@parboxrestore\tiny\raggedright #2}}
\catcode`\@=12
\newcommand{\intersect}{\cap}
\newcommand{\union}{\cup}
\newcommand{\CC}{{\mathcal C}}
\newcommand{\til}{\widetilde}
\newcommand{\subtrack}{<}
\newcommand{\ML}{\mathcal{ML}}

\long\def\realfig#1#2{
\begin{figure}[htbp]
\centerline{\psfig{file=#1.ps}}
\caption[#1]{#2}
\label{#1}
\end{figure}}

\long\def\state#1#2{
\medskip\par\noindent
{\bf #1} 
{\it #2}
\par\medskip
}

\author{ Howard Masur and Yair Minsky}
\date{\today}
\begin{document}

\title{Quasiconvexity in the curve complex}
\maketitle

\section{Introduction}
Let $S$ be a surface of genus  $g$ with $n$ punctures, and assume
$3g-3+n>1$.  Associated 
to $S$ is an object $\CC(S)$ called the {\em complex of curves}, whose
vertices are homotopy classes of nontrivial, nonperipheral simple closed
curves. A $k$-simplex of $\CC(S)$ is a collection of $k+1$ disjoint
nonperipheral homotopically distinct simple closed curves. The
dimension of the complex is $3g-4+n$. We are interested in the {\em
  $1$-skeleton} $\CC_1(S)$ which is endowed with a metric
$d(\cdot,\cdot)$ which assigns  length $1$ to every edge.

In \cite{masur-minsky:complex1} we proved that $\CC_1(S)$, with this metric, is
$\delta$-hyperbolic in the sense of Gromov and Cannon. This fact has
applications to the study of mapping class groups and hyperbolic
3-manifolds
\cite{masur-minsky:complex2,minsky:ELCI,brock-canary-minsky:ELCII}. 
See also Bowditch \cite{bowditch:complex} for a new and more concrete proof of the
hyperbolicity theorem.

In the study of Heegaard splittings it is of interest to consider the
set of essential curves in the boundary of a handlebody which bound 
disks -- see Casson-Gordon
\cite{casson-gordon:reducing}, Hempel \cite{hempel:complex}, Johannson
\cite{johannson:3mfds} and Schleimer \cite{schleimer:distance2}. 
In slightly more generality, for a compact, orientable
3-manifold $M$ with boundary component $S$ we can define the set of
(homotopy classes of) boundaries of essential disks: 
$$
\Delta(M,S) = \{[\partial D]: (D,\partial D)\subset (M,S) \ 
\text{is an essential disk}\} \subset \CC(S).
$$
The purpose of this note is to prove

\begin{theorem}
\label{theorem:handlebody}
$\Delta(M,S)$ is a $K$-quasiconvex subset of $\CC(S)$, where $K$
depends only on the genus of $S$.
\end{theorem}

Recall that a subset $Y$ of a geodesic metric space $(X,d)$ is {\em
$K$-quasiconvex} if for any pair of points $y_1,y_2\in Y$ any geodesic
in $X$ joining them stays in a $K$-neighborhood of $Y$.  This theorem can be
  viewed as something of an analogy to a result on
quasiconvexity of sublevel sets in $\CC(S)$ for certain length
functions arising from Kleinian representations, which plays a central role 
in \cite{minsky:boundgeom,minsky:ELCI}.

Theorem \ref{theorem:handlebody}
will be a special case of more general theorems about quasi-convex
  subsets of $\CC(S)$ obtained by the combinatorial processes of 
{\em curve replacement} and
{\em train-track nesting}.

Curve replacement is the standard process of simplifying intersections
with a closed curve in a surface by surgery, which is usually used to
obtain upper bounds on distance in $\CC(S)$. We will treat this with a
bit of care in Section \ref{replacements}, because we need to consider
{\em well-nested curve replacements}. In Section \ref{prove curve
  replacements} we will prove:  

\begin{theorem}\label{Curve replacements}
There exists $K=K(S)$ such that the vertices of a 
well-nested curve replacement sequence form a 
$K$-quasiconvex set in $\CC(S)$, which moreover is Hausdorff distance
$K$ from a geodesic in $\CC(S)$. 
\end{theorem} 

In order to prove this theorem we will consider {\em nested train
tracks} in Section \ref{train tracks}, and prove the following
statement (stated more precisely later on):

\begin{theorem}\label{tracks quasiconvex}
The vertices of a nested train-track sequence with bounded steps form
a $K$-quasiconvex set in $\CC(S)$, which moreover is Hausdorff distance
$K$ from a geodesic in $\CC(S)$.  
\end{theorem}

These two theorems, together with Proposition \ref{prop:nested} which
produces curve-replacement sequences lying in $\Delta(M,S)$, will be enough
for the proof of  Theorem \ref{theorem:handlebody}, which will appear at the end
of  Section \ref{prove curve replacements}.

\section{Curve Replacement}\label{replacements}
A closed curve in $S$ is {\em essential} if it is homotopically
nontrivial and nonperipheral (not homotopic to a puncture of $S$).
Two essential curves $a$ and $b$ are in 
{\em minimal position} if they intersect transversely and the number
of intersection points 
$|a\intersect b|$ is minimal over all representatives of their
respective homotopy classes. This number is the geometric intersection
number of the homotopy classes, written $i([a],[b])$ or just
$i(a,b)$. 

Consider two simple curves $a$ and $b$ in minimal position, and let
$J\subset a$ be an interval with $J\intersect b \ne \emptyset$. 
For any $x,y\in J$ let $[x,y]$ denote the subinterval
of $J$ with endpoints $x$ and $y$, and similarly $(x,y)$,  $(x,y]$
etc. 

A {\em curve replacement} with respect to $(a,b,J)$ is the 
construction of a new curve $a_1$ and a new interval $J_1\subset J$ in
one of the following ways: 
\begin{enumerate}
\item Let $w$ be an interval of $b$ with $int(w)\intersect a  =
  \emptyset$, and endpoints $p,q\in int(J)$. 

  If $w$ is incident to opposite sides of $a$ at its two endpoints, 
  Let $J_1=[p',q']$ where $p'$ is
  slightly to the right of $p$ and $q'$ is slightly to the right of $q$,
  so that $J_1\intersect b = (p,q]\intersect b$. Let $w'$ be a slight
  perturbation of $w$ in a regular neighborhood so that the endpoints
  of $w'$ are $p',q'$ and $w'$ and $w$ are disjoint.
  Define $a_1$ to be the composition $w' * J_1$. 
\realfig{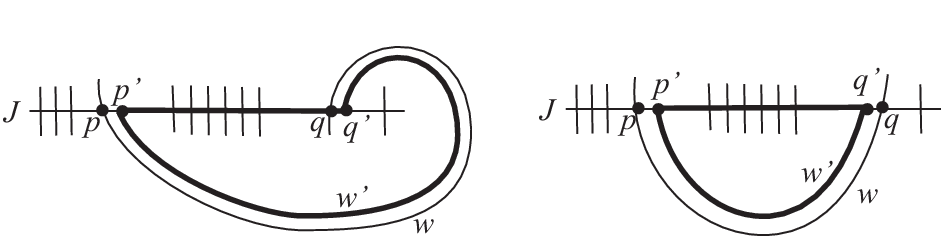}{The two curve replacements using a single arc $w$. 
In each case $a_1$ is thickened.}

  If $w$ is incident to the same side of $a$ at both endpoints, we
  call $w$ a {\em wave}. In this case, let us also assume that $w$ is
   {\em innermost}, i.e. that there is no other wave $w_2$ with
   endpoints inside $[p,q]$ incident to $w$ on the same side of $J$
   such that $w$ and $w_2$ are homotopic with 
   endpoints on $J$.   In this case let $J_1=[p',q']$ where
   $p',q'\in(p,q)$ and $J_1\intersect b = (p,q)\intersect b$. Again
  move $w$ in a regular neighborhood to obtain $w'$ with endpoints
  $p',q'$, and let $a_1 = w' * J_1$.

\item Let $w_1$ and $w_2$ be two waves incident to  $J$ on opposite
  sides of $a$, with endpoints $p_1,q_1$ and $p_2,q_2$ respectively.
  Assume that $p_1<q_1$ and $p_2<q_2$ after fixing some orientation on
  $J$. Suppose either that $q_1=p_2$ or $q_1=q_2$.
   Move each of $w_1$ and $w_2$ slightly inward, as above;  in the
case $q_1=q_2$ do this so that $q_1'=q_2'$. In either case we
obtain $w'_i$ with endpoints 
   $p'_i<q'_i$. Now let $a_1 = w'_1 * [p'_1,p'_2] * w'_2 *
   [q'_2,q'_1]$. Let $p' = \min(p'_1,p'_2)$.  In the case that $q_1=p_2$ 
 set  $q'=q'_2$, and in the case of $q_1=q_2$ set 
$q'=\max(p_1',p_2')$. Finally set  
   $J_1 = [p',q']$.

  We call this a {\em double wave} curve replacement. 

\realfig{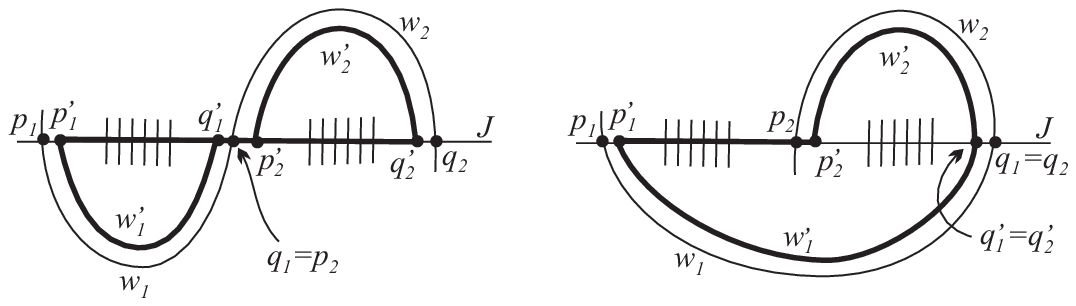}{The two types of double wave replacements. Note that
in the $q_1=p_2$ case $a_1$ is not embedded, but can be perturbed to be.}

\end{enumerate}
Note that in case (1), and case (2) when $q_1=q_2$, $a_1$ is simple,
and in case (2) when $q_1=p_2$ it is homotopic 
to a simple curve (the arcs $[p_1,p_2]$ and $[q_1,q_2]$ 
overlap  but after moving them  slightly to opposite sides of $a$ we
get a simple curve).
In all cases $a_1$ is homotopically nontrivial,
which follows from the assumption that $a$ and $b$ are in minimal
position. 

In the non-wave and double wave case it is easy to check that $a_1$ is
also nonperipheral. If $S$ is closed then of course $a_1$ is
always nonperipheral.

We note the special case of a double-wave replacement with $p_1=p_2$
and $q_1=q_2$.  In this case  $i(a,b)=2$ and $a_1$ is homotopic to $b$. 

One can also check that $a_1$ and $b$ are in minimal position -- here
in the single wave case we must use the condition that the wave is 
innermost. 

Note that all intersections of $a_1$ with $b$ lie in $J_1$. 
For the curve replacements in case (1) and the double wave replacement
when $q_1=q_2$, we have in particular
\begin{equation}\label{int number dec 1}
i(a_1,b) \le |J_1\intersect b| < |J\intersect b|.
\end{equation}
For the double wave replacement when $q_1=p_2$ the counting is
slightly different and we get:
\begin{equation}\label{int number dec 2}
  i(a_1,b) \le |J_1\intersect b|+1 < |J\intersect b|.
\end{equation}

We also define, for completeness, a curve replacement in the case that
$J\intersect b = \emptyset$ to be $a_1=b$, and $J_1=\emptyset$. If
$J\intersect b$ contains 
exactly one point then we are in the non-wave case of part 1 with
$p=q$, and the replacement is again $a_1=b$, $J_1=\emptyset$.

\subsection*{Nested curve replacement}
Given $a$ and $b$ in minimal position, a {\em nested curve replacement
sequence} is a sequence $\{(a_i,J_i)\}$ of curves $a=a_0,a_1,a_2,\ldots,a_n$ and
segments $a\supset J_0\supset J_1 \supset \cdots \supset J_n$ such that
$J_0$ contains all the points of $a\intersect b$, and
$a_{i+1},J_{i+1}$ are obtained by a curve replacement from 
$(a,b, J_i)$.

\begin{proposition}
\label{prop:nested}
Given any $a,b$ in minimal position, and an interval $J_0\subset a$
containing $a\intersect b$,
there exists a nested curve replacement sequence  $\{(a_i,J_i)\}$
such that
\begin{itemize}
\item  $a_i$ is nonperipheral for all $i$.    
\item If $S$ is a boundary component of a  compact 3-manifold $M$ and
$a$ and $b$ are boundaries of compressing discs, then the $a_i$
can be chosen to be boundaries of compressing discs. 
\item The sequence terminates with $a_n$ homotopic to $b$.
\end{itemize}
\end{proposition}
\begin{proof} The sequence is chosen inductively.   
To satisfy the first statement, 
at each stage we choose a curve replacement configuration as in the
definition, noting that if there are no non-wave curve replacements
and every single wave replacement yields a peripheral curve, then
there must be a double wave replacement, which will yield a
non-peripheral curve. The intervals $J_i$ are automatically nested by
the definition. The process terminates, with $a_n$ homotopic to $b$,
because the numbers 
$|J_i\intersect b|$ are strictly decreasing by (\ref{int number dec
1}) and (\ref{int number dec 2}).

To prove the second statement, note first that since $S$ is closed
the non-peripheral property is automatic. We will inductively
construct each $a_i$ using a wave curve replacement. 
We recall from the definition of curve replacements and induction, 
that for $i>0$, $a_i$ intersects $b$ 
transversely and $a_i\intersect b\subset J_i$ (for $i=0$
this is true by hypothesis). 
Let $B$ be a properly embedded disk in $M$ with boundary $b$, and 
suppose by induction that there exists a properly embedded disk $A_i$
with boundary $a_i$. If $i(a_i,b)=0$ we are done, so suppose
$i(a_i,b)>0$. 

We may assume that $A_i$ and $B$ intersect transversely.
Their intersection locus can be assumed to consist of properly
embedded arcs, since the closed-curve components can be removed by an
exchange that does not affect $a_i$ and $b$. Now let $e$ be an
innermost intersection arc on $B$, i.e. an arc that, together with
an arc $w$ on $b$ bounds a disk $E$ in $B$ whose interior is disjoint
from $A_i$. Thus $int(w)$ is disjoint from $a_i$ and hence from
$J_i$. 
The endpoints $p,q$ of $w$ lie on $a_i$, and hence on $J_i$. 
The arc $[p,q]\subset J_i$ together with $e$ bounds a subdisk of $F$ of
$A_i$. The disk $E\union F$ has boundary $[p,q]*w$ which, after a
slight isotopy, is the curve replacement $a_{i+1}$ obtained from $(a,b,J_i)$
using the wave $w$. 
\end{proof}

\subsection*{Distance in $\CC(S)$}

Using curve replacements one can bound distance in $\CC(S)$ in terms
of intersection number. One such bound is:
\begin{equation}
  \label{dist bound}
  d(\alpha,\beta) \le i(\alpha,\beta)+1
\end{equation}
(see Lemma 1.1 of Bowditch \cite{bowditch:complex}. For large
intersection numbers a better bound is logarithmic in
$i(\alpha,\beta)$ -- see also Hempel \cite{hempel:complex} and
Masur-Minsky \cite{masur-minsky:complex1}.)

If $\{(a_i,J_i)\}$ is a curve replacement sequence, we note that 
$| a_{i+1} \intersect J_i| \le 2$. It follows (taking a little care
when there are two
double-wave replacements in a row) that
\begin{equation}
  \label{eq:succesive intersection}
  i(a_i,a_{i+1}) \le 4.
\end{equation}
From this and (\ref{dist bound}) we obtain
\begin{equation}
  \label{eq:successive distance}
  d(a_i,a_{i+1}) \le 5.
\end{equation}

We note that (\ref{eq:successive distance}) is probably not sharp, even
in the case of double-wave 
replacements, but it will suffice for our purposes. 
In the case of one-wave replacements (the only case relevant to 
Theorem \ref{theorem:handlebody}), 
we always have $i(a_i,a_{i+1})=0$ and $d(a_i,a_{i+1})\le 1$.

\section{Train tracks}
\label{train tracks}

A train track on a surface $S$ is an embedded $1$-complex $\tau$
satisfying the following properties (see Penner-Harer
\cite{penner-harer} for more details). Each edge (called a branch) is a
smooth path with well-defined tangent vectors at the endpoints, and at
any vertex, (called a switch) the incident edges are mutually tangent.
The tangent vector at the switch pointing toward the interior of an
edge can have two possible directions, dividing the ends of edges into
``outgoing" and ``incoming" directions, neither of which can be empty.
For the tracks that we consider there will be a single switch. Each
component of $S\setminus \tau$ has negative generalized Euler
characteristic which is the usual Euler characteristic minus $1/2$ for
every outward pointing cusp (internal angle $0$).  This means that we
exclude annuli, once-punctured discs with smooth boundary, or
unpunctured discs with $0$,  $1$, or $2$  cusps.  The last are called bigons. 

A {\em train route} is a nondegenerate smooth path in $\tau$. It
traverses a switch only by passing from incoming to outgoing edge or
vice versa. A {\em transverse measure} on $\tau$ is a non-negative
function $\mu$ on branches satisfying the switch condition that for
any switch the sums of $\mu$ over outgoing edges equals the sum over
incoming. A closed train route induces counting measure on $\tau$.

A train track is {\em recurrent} if every branch is contained in a
closed train route. A curve $\beta$ is {\em carried} on $\tau$ if
there is a homotopy of $S$ taking $\beta$ to a set of train
routes. Then $\beta$ induces a measure on $\tau$ which uniquely
determines it.

 A train track $\sigma$ is {\em carried} by $\tau$, written
$\sigma\prec \tau$,  if there is a homotopy of the surface such that
every train route of $\sigma$ is taken to a train to a train route of
$\tau$.  We also say that $\sigma$ is {\em nested} in $\tau$.
We write $\sigma<\tau$ if $\sigma$ is a {\em subtrack} of
$\tau$; that is, $\sigma$ is a train track which is a subset of
$\tau$.

For a recurrent train track $\tau$ let $P(\tau)$ denote the polyhedron
of transverse measures supported on $\tau$. We note that $P(\tau)$ is preserved
by scaling so it is a cone on a compact polyhedron in a projective
space.  We can consider all of these polyhedra as subsets of the
measured lamination space $\ML(S)$, and then
$\sigma\prec \tau$ is equivalent to $P(\sigma)\subset P(\tau)$,
and $\sigma<\tau$ is equivalent to $P(\sigma)$ being a face of $P(\tau)$. 
By $ int(P(\tau))$ we will mean the set of measures on $\tau$
which are positive on every branch (recurrence implies that this is nonempty).

A {\em vertex} of the projectivization of $P(\tau)$ is an extreme
point, and corresponds to a line in $\bdry P(\tau)$. By abuse of
notation we will let ``vertex'' denote any non-zero point on this
line. After scaling, a vertex can always be  realized by the
counting measure on a single simple closed curve (see \cite[\S
4.1]{masur-minsky:complex1}), which we call a {\em vertex cycle.}

Because the set $vert(\tau)$ of vertex cycles of $\tau$ is finite
and there are finitely many
homeomorphism classes of train tracks in $S$, 
there is a constant $B$ such that 
$$diam_{\CC(S)} \{vert(\tau)\}\leq B.$$

For tracks $\tau$ and $\sigma$, define
\begin{displaymath}
d_T(\tau,\sigma)=\min_{\substack{\beta\in vert(\tau)\\\alpha\in vert(\sigma)}}
d(\beta,\alpha)
\end{displaymath} 
Note that this is not a distance function on the set of tracks. 

Our goal in this section is to prove Theorem \ref{tracks quasiconvex},
which we restate here more precisely:
\state{Theorem \ref{tracks quasiconvex}.}{
Given $M>0$ there exists $K$ such that if  $\cdots\prec\tau_n\prec
 \tau_{n-1}\prec \cdots\prec\tau_1\prec\tau_0$ is a sequence of 
nested train tracks such that $d_T(\tau_{i+1},\tau_i)\leq M$  
then the set of vertices of the $\tau_i$ is $K$-quasiconvex in $\CC(S)$.
}

\medskip

A track is {\em large} if all complementary components of the track are
disks or once-punctured disks. We have:
\begin{lemma}\label{big is large}
If a train track $\tau$ is not large then for any
$\alpha,\beta$ carried by $\tau$,
$d(\alpha,\beta)\le 2$.
\end{lemma}
\begin{proof}
Since $\tau$ is not large, there is an essential simple closed curve $\gamma$
which misses $\tau$, and hence both $\alpha$ and $\beta$; hence
$d(\alpha,\beta)\le 2$.  
\end{proof}

This lemma will mean that we will be able to restrict our attention to
large tracks.  Let $\sigma$ be a large track. A {\em diagonal
  extension} of $\sigma$ is a track $\kappa$ such that $\sigma<\kappa$
and every branch of $\kappa\setminus \sigma$ is a diagonal of
$\sigma$.  This means that it is an edge that terminates in corners of
a complementary region of $\sigma$.

Let $E(\sigma)$ denote  the (finite) set of all recurrent tracks which are diagonal extensions of $\sigma$. Let 
\begin{displaymath}
PE(\sigma)=\bigcup_{\kappa\in E(\sigma)} P(\kappa)
\end{displaymath}
Let $int(PE(\sigma))$ denote the set of measures $\mu\in  PE(\sigma)$
which are positive on every branch of $\sigma$.
Next for $\tau$ a large track, let 
\begin{displaymath}
N(\tau)=\bigcup_{\substack{\sigma<\tau\\ \text{$\sigma$ large}}} E(\sigma),
\end{displaymath}
and define
$$
PN(\tau)=\cup_{\kappa\in N(\tau)}P(\kappa).
$$
This can be thought of as a neighborhood of $P(\tau)$. Let
\begin{displaymath}
int(PN(\tau))=\bigcup_{\substack{\sigma<\tau\\ \text{$\sigma$ large}}}
int(PE(\sigma)).
\end{displaymath} 

We will need the following  lemmas proved in  \cite{masur-minsky:complex1}. 
\begin{lemma}\label{lemma:deeply:nested} 
 There is a constant $D$ such that if 
$\sigma\prec\tau$ are a pair of large recurrent tracks and $d_T(\tau,\sigma)\geq D$, then 
$PN(\sigma)\subset int(PN(\tau))$
\end{lemma}
This lemma is stated  in \cite{masur-minsky:complex1} for {\em generic tracks}, i.e. those
with trivalent switches, but any track can be perturbed to a generic
track without changing its vertex set or polyhedron of measures
(by ``combing'' -- see Penner-Harer \cite[\S1.4]{penner-harer}), and
it is easy to see that this can be done simultaneously for $\sigma$
and $\tau$.

\begin{lemma}\label{lemma:extension}
If $\sigma$ is a large track, $\gamma\in int(PE(\sigma))$ and
$d(\gamma,\beta)\leq 1$, then $\beta\in PE(\sigma)$.  
\end{lemma}
In fact this implies 
\begin{lemma}\label{lemma:PN nest}
If $\sigma$ is a large track then
$$
\NN_1(int(PN(\sigma))) \subset PN(\sigma)
$$
\end{lemma}
where $\NN_k(int(PN(\sigma)))$ denotes the $k$-neighborhood in $\CC(S)$ of
the vertices of $\CC(S)$ contained in $int(PN(\sigma))$.

\medskip

We will also need the following lemma. 
\begin{lemma}
\label{lemma:disjoint:interiors}
Let $\tau$ be a large recurrent track and $v$ a vertex of $\tau$. 
Then $v\notin int(PN(\tau))$.
\end{lemma}
\begin{proof}
First we claim  that $dim (P(\tau))\geq 2$ for any large recurrent track $\tau$.  
Let $b$ be the number of branches,  $s$ the number of switches and  $f$
the number of complementary discs. Then  
\begin{displaymath}\chi(S)=f-b+s
\end{displaymath} so 
\begin{displaymath}
dim (P(\tau))\geq b-s=f-\chi(S),
\end{displaymath}
since $b-s$ is a lower bound for the dimension of the solution set to
the switch conditions
and the recurrence of $\tau$ gives us at least one
solution with all positive weights.
The latter quantity is at least $2$ except in the case that $S$ is a
punctured torus, where $\chi(S)=-1$ and $f=0$.  However in this case
it is easily checked directly that  $dim (P(\tau))=2$, so the claim is
proven.

Suppose now that the lemma is false. By definition then, there must be
some large subtrack $\sigma<\tau$, and track $\kappa\in E(\sigma)$
such that a representative  $\alpha_1$ of $v$
is carried by $\kappa$ and assigns positive weights to the
branches of $\sigma$. 
We may assume that $\alpha_1$ assigns positive
weights to every branch of $\kappa$ as well -- deleting branches of
$\kappa$ if necessary. Let $\alpha_2$ be a representative of $v$
carried on $\tau$. We may assume that $\alpha_2$ puts positive measure
on every branch of $\tau\setminus\sigma$, otherwise we may remove
those branches, obtaining a smaller track containing $\sigma$ for
which $v$ is a vertex. We now wish to arrive at a contradiction.   

By an {\em admissible deformation} of an edge $e$ of $\kappa\setminus
\sigma$, we mean a homotopy, preserving endpoints and tangency to
$\sigma$, to a path $e'$ such that $\sigma\union e'$ is still a train
track. In particular $e'$ is allowed to have segments tangent to $\sigma$.

We say an edge $e$ of $\kappa\setminus\sigma$
intersects an edge $f$ of $\tau\setminus\sigma$ {\em inessentially} if $e$ can be
deformed admissibly to an edge $e'$ which is either
disjoint from the interior of $f$, or traverses a train route in
$\sigma\union f$. Suppose that all intersections of $\tau\setminus
\sigma$ with $\kappa\setminus\sigma$ are inessential. Then $\kappa$
can be deformed to be carried in a track $\hat\tau$ containing $\sigma$,
on which $v$ puts
a measure positive on every edge, and such that $\tau$ and $\hat\tau$
are subtracks of a common train track $\omega$. 

Suppose first that $\hat\tau\subtrack\tau$.
Since $\sigma$ is large so is $\hat\tau$, but then by the first claim
$dim(P(\hat\tau))\ge 2$, which means the projectivized polyhedron has
dimension at least 1, and any extreme point could not be in its interior.
Thus $v\in int(P(\hat\tau))$ means that $v$ is not extreme in
$P(\hat\tau)$, and hence not in $P(\tau)$.

Next suppose that $\hat\tau$ is not a subtrack of $\tau$. In this case
$v$ is represented by two distinct measures on $\omega$, which
contradicts the injectivity of the map $P(\omega)\to \ML(S)$.

We are left with the possibility that an edge $b_\tau$ of
$\tau\setminus\sigma$ has essential intersection with $b_\kappa$ in
$\kappa\setminus\sigma$. Both edges lie in some
complementary disc or once punctured disc $R$  of $S\setminus\sigma$,
and we note that $b_\kappa$ is a diagonal but $b_\tau$ does not need
to be. The representative $\alpha_2$ of $v$ 
assigns positive weight to $b_\tau$, and the
representative $\alpha_1$ assigns positive weight to
$b_\kappa$. 

We may assume that $b_\tau$ and $b_\kappa$ have a unique intersection
point $x$ if $R$ is a disk, or possibly two intersection points if $R$
is a punctured disk. Let $x$ be one of the points in the latter case. 
We may also assume, possibly after admissible deformations of
$\kappa$,  that there are no inessential intersections of
$\kappa$ and $\tau$, and that every edge of $\kappa$ that is
admissibly deformable into $\tau$ already lies in $\tau$. 

Let $\til\tau$ and $\til\kappa$ be the lifts of the tracks $\tau$
and $\kappa$ to the
universal cover.  Each contains the lift $\til\sigma$ of $\sigma$.
Let $\til x\in \til R$ be lifts of $x$ and $R$.  Let
$\til\alpha_1,\til\alpha_2$ be lifts of $\alpha_1,\alpha_2$ that pass
through $\til x$, and let $\til b_\tau,\til b_\kappa$ be the lifts
of $b_\tau$ and $b_\kappa$ that intersect at $\til x$. 
Since $\alpha_1\sim\alpha_2$, $\til\alpha_1$
and $\til\alpha_2$ must intersect at some other point $\til y$ so
that the segments of $\til\alpha_1$ and $\til\alpha_2$ from
$\til x$ to $\til y$ bound a disc $D$.  The point $\til y$ may
or may not be in the interior of a complementary domain of
$\til\sigma$.

We will now build a track $\til\sigma_1$ that contains
$\til\sigma$.  To form $\til\sigma_1$, add to $\til\sigma$ any
branch of $\til\alpha_1$ and $\til\alpha_2$ that is entirely
contained in $\bdry D$. Since $\bdry D$ is embedded none of
these branches cross and we obtain a train track. 
Notice that we do not introduce
any bigons, since neither $\tau$ nor $\kappa$ have any, and we have
already deformed $\kappa$ so that no edge of $\kappa\setminus\tau$ is admissibly
deformable into $\tau$.  Thus any
complementary disc of $\til\sigma_1$ still has 
negative generalized Euler characteristic.

Next consider $\til R\cap D$.  It is bounded by one or more edges of
$\til\sigma$, a subsegment of $\til b_\tau$ that goes from a
switch $\til P_1$ to $\til x$ and a subsegment of $\til
b_\kappa$ from $\til x$ to a switch $\til P_2$.  Replace this
latter pair of segments with a smooth path joining $\til P_1$ to
$\til P_2$ in the same homotopy class rel endpoints (and maintaining
tangency at the switches).  The assumption
that the intersection is essential implies that the subdisc $\tilde R'$ of $
D$ bounded by this smooth path and the edges of $\til\sigma$ is not
a monogon. If $\tilde R'$  is a not a bigon, we add this edge to $\til
\sigma_1$ and replace 
$D$ with $D'=(D\setminus \til R)\cup \til R'$.
Note that $\til\sigma_1 $ still
has no bigons. If $\tilde R'$  is a bigon,
the path can be admissibly deformed into
$\til\sigma$ and we do not add it to $\til\sigma_1$.
Now let $D'=D\setminus \til R$.

We perform the same construction if
$\til y$ is in the interior of a complementary domain as well,
noting that $\til y$ is also an essential intersection point. 
The new disc $D''$ has smooth boundary except possibly for
one cusp at the point $\til y$ (if $\til y$ was not an interior
point). Thus its generalized Euler
characteristic is positive.  However it is a union of complementary
discs of $\til\sigma_1$, each of which has negative generalized Euler
characteristic.  Since generalized Euler characteristic adds, we have a
contradiction.
\end{proof}

We are now ready to give the proof of the quasiconvexity theorem for
nested sequences of tracks. 

\begin{proof}[Proof of Theorem \ref{tracks quasiconvex}]
Let $\{\tau_n\}$ be our nested sequence, so that
$\tau_{i+1}\prec\tau_i$
and $d_T(\tau_i,\tau_{i+1})\le M$. 
We may assume that all $\tau_i$ are large. For if some $\tau_i$ is not
large so are all $\tau_j$ for $j\ge i$, and since all are carried in
$\tau_i$, Lemma \ref{big is large} implies that $d_T(\tau_j,\tau_k) <
3$ for all $j,k\ge i$, and hence we may ignore these tracks. 
Using the number $D$  given by Lemma~\ref{lemma:deeply:nested}, the
condition $d_T(\tau_i,\tau_{i+1})\le M$
implies that  we may inductively find a subsequence
$\tau_{i_j}$ (with $i_0=0$) such that  
\begin{equation}\label{tauij quasi}
D \leq d_T(\tau_{i_j},\tau_{i_{j+1}}) < D+M
\end{equation}
and for any $\tau_n$ we have
\begin{equation}\label{small}
d_T(\tau_{i_j},\tau_n) < D
\end{equation}
for some $i_j$.

Lemma~\ref{lemma:deeply:nested} then says that
$$PN(\tau_{i_{j+1}})\subset int(PN(\tau_{i_j})),$$  and
Lemma~\ref{lemma:PN nest} says that
$$\NN_1(int(PN(\tau_{i_{j}}))) \subset PN(\tau_{i_j}).$$
Combining these and applying induction we find that
$$
\NN_{k-1}(PN(\tau_{i_{j+k}})) \subset int(PN(\tau_{i_j})).
$$
Applying this to the vertices, and using Lemma~\ref{lemma:disjoint:interiors} 
 we find
\begin{equation}\label{distance:bound}
d_T(\tau_{i_j},\tau_{i_k})\geq |k-j|
\end{equation}
for any $j,k$ in the subsequence.

Inequality (\ref{tauij quasi}), together with the bound $diam
\{vert(\tau)\}\leq B$, implies that any sequence of vertices of the
$\tau_{i_j}$ forms a $K$-quasigeodesic with $K=K(D,M,B)$. By
$\delta$-hyperbolicity of $\CC(S)$, any geodesic joining the vertices
of this sequence stays within some fixed distance $M'(M,B,D,\delta)$
of it. Thus $\bigcup vert(\tau_{i_j})$ is $M'$-quasi-convex.

Now by (\ref{small}), $vert(\tau_i)$ for any $i$ is within $D+B$ of some
$vert(\tau_{i_j})$. Again by $\delta$-hyperbolicity, there exists
$M''=M''(B,D,\delta)$ such that any
geodesic joining vertices of the $\tau_i$ stays within distance $M''$
of a geodesic joining some pair of vertices of the $\tau_{i_j}$.

Combining
these two facts we see that $\bigcup vert(\tau_i)$ is a $M'+M''$
quasiconvex set, and moreover that it is Hausdorff distance $M'+M''$
from a geodesic connecting vertices of $\tau_{i_0}$ and the last
$\tau_{i_j}$. 
\end{proof}

\section{Quasiconvexity of curve replacements}
\label{prove curve replacements}
In this section we will restate Theorem \ref{Curve replacements} more
precisely, and 
apply Theorem~\ref{tracks quasiconvex} to prove it. 
In order to do this, given a nested
curve replacement sequence $(a_i,J_i)$ we 
will construct a sequence of nested train tracks $\tau_i$ carrying
$b$ such the vertices of $\tau_i$ are a bounded $\CC(S)$-distance from
$[a_i]$.

\subsection*{One-switch train tracks}
Let $a$ and $b$ be two essential curves intersecting minimally and let
$J$ be an interval in $a$ such that $J\intersect b$ is non-empty.
We wish to construct a train track $\tau_J$ 
closely associated to the ``first return map'' of $b$ to $J$.
The quotient of $J\union b$ obtained by squeezing $J$ to a point $P$ is
a 1-complex, and we may impose a switch structure on $P$
by declaring that tangent vectors at $P$ pointing to interiors are all
outgoing for edges leaving $a$ on one side and incoming for edges
leaving $a$ on the other.  This makes the quotient,
which we call $\tau^0_J$,  a
``bigon track'': it satisfies all the
conditions for being a train track except for the possibility of bigon
complementary regions. The other disallowed types of complementary
regions fail to occur because $a$ and $b$ intersect minimimally.

If we identify the opposite sides of any bigon in $\tau^0_J$ we obtain
a new 1-complex $\tau_J$, still with a switch structure. This will be
a train track provided that the identification can be performed as a
homotopy of the surface $S$, and in this case we say that $J$ is a
{\em good interval}. The following lemmas will allow us to use this
construction effectively:

\begin{lemma}\label{good intervals}
Let $a$ and $b$ be essential curves in $S$ intersecting
minimally. 
\begin{enumerate}
\item There exists a good interval  $J\subset a$ containing all points of
$a\intersect b$.
\item Any subinterval of a good interval is good. 
\end{enumerate}
\end{lemma}
\begin{proof}
A component $Q$ of $S\setminus a\union b$ is a {\em rectangle} if it
is a disk and its boundary, traversed from the inside, is a path
meeting exactly four intersection points of $a\intersect b $. By an Euler
characteristic argument, there must be at least one region $Q$ which
is not a rectangle. Let 
$I\subset a$ be an open interval between two successive points $x_1,x_2$ of
$a\intersect b$, which lies on the boundary of a non-rectangle $Q$.
We claim that $J=a\setminus I$ is a good interval. 

 To see this, squeeze $J$ to a point and let $\Omega$ be the
 complementary region of the resulting $\tau_J^0$ 
which contains $Q$. For each
$i=1,2$ the arc of $b$ that passes through $x_i$ is on the boundary
 of $\Omega$ and traverses the switch 
$P$ smoothly.  If $\Omega$ were a bigon, its boundary would consist of
 two arcs of $b$ passing smoothly through $x_1$ and $x_2$ and meeting
 in cusps at both ends.  The preimage of these cusps would consist of 
two arcs of $b\setminus a$, and this would force $Q$ to be a rectangle, a
 contradiction. 
 Furthermore, because $\tau_J^0$ has a single switch $P$, 
 the arcs of $b$ through $x_i$ are the only smooth arcs
 through $P$ on the boundary of complementary regions. As they are on
 the boundary of $\Omega$ which is not a bigon, every  bigon of
 $\tau^0_J$ meets the switch only at its 
cusps, and hence has exactly two branches of $\tau^0_J$ in its
boundary. Consider a connected component of the union of bigons in the
complement of the switch. It must be a disk (with punctured boundary)
with two branches in its boundary -- the only other way of gluing
such bigons together gives a sphere minus two points, which is impossible. Hence the
collapse of bigons can be performed on each connected component, yielding the
desired track $\tau_J$. 

For part (2): Suppose $J\subset a$ is good, 
let $J'\subset J$ be an interval which intersects
$b$, and suppose $B$ is a bigon complementary region of
$\tau^0_{J'}$. Every arc of $b\setminus J'$ is a train route for
$\tau^0_J$, so clearly $\tau^0_{J'}\prec \tau^0_J$. Thus after a
homotopy of $S$ we can assume the branches of $\tau^0_{J'}$ travel
along branches of $\tau^0_J$, and $B$ is taken by this homotopy to a
collection of complementary regions of $\tau^0_J$, attached smoothly
along arcs of $\tau^0_J$. 
The generalized Euler characteristic of any complementary region 
is nonpositive and  zero if and only if it is a bigon or an annulus.
Since $B$ is simply connected and the
characteristic of a union of complementary regions is additive, 
these  regions must necessarily be bigons. Hence the bigon
collapse $\tau^0_J \to \tau_J$ also collapses $B$. Since this collapse
can be performed by a homotopy of $S$, this same homotopy collapses
all the bigons of $\tau^0_{J'}$, so $J'$ is good. 
\end{proof}

With this in mind, let us say that $\{(a_i,J_i)\}$  is a {\em
  well-nested curve replacement sequence} if it is a nested curve
  replacement sequence and $J_0\subset a$ is a good interval. We can
  now restate Theorem \ref{Curve replacements}.
\state{Theorem \ref{Curve replacements}.}{
There exists $K=K(S)$ such that, if $a$ and $b$ are in minimal
  position and $\{(a_i,J_i)\}$ is a well-nested curve replacement
  sequence, then the vertices $[a_i]$ form a
$K$-quasiconvex set in $\CC(S)$ which moreover are Hausdorff distance $K$
 from a geodesic in $C(S)$. }

\begin{proof}
Now let $J_0\subset a$ be a good interval which contains $a\intersect
b$, 
and let $\{(a_i,J_i)\}$ be a curve replacement sequence, terminating
with $[a_n]=[b]$.
Using Lemma \ref{good intervals}, each $J_i$ is good, and so we
obtain a sequence of tracks $\tau_i = \tau_{J_i}$. 
As in the proof of Lemma \ref{good intervals}, we have: 
\begin{lemma} For each $i$, 
$\tau_{i+1}\prec\tau_i$, and furthermore each $\tau_i$ carries $b$.
\end{lemma}

Next we observe that 
\begin{equation}\label{alpha vertex bound}
d(a_i,vert(\tau_i)) \le 5
\end{equation}             

To prove (\ref{alpha vertex bound}), we observe that $a_i$ is
composed of one or two arcs of
$b$ and one or two arcs along $J_i$. Hence after after isotopy it
intersects $\tau_i$ at most 
twice at the switch $P$. Each vertex cycle of $\tau_i$ passes through
the switch at most twice.
Hence $i(a_i,v)\le 4$ for any vertex $v$ of $\tau_i$, and the distance
bound follows from (\ref{dist bound}).

Since, by (\ref{eq:successive distance}), $d(a_i,a_{i+1}) \le 5$, we may conclude by the
triangle inequality that
$$d_T(\tau_i,\tau_{i+1}) \le 15.$$ 
(As in  Section \ref{replacements}, we note that this is probably not
sharp. In the case of one-wave replacements the bound obtained is 7.)

Thus the hypotheses of 
Theorem \ref{tracks quasiconvex} are satisfied, so that
$\{vert(\tau_i)\}$ form a quasiconvex set. Inequality (\ref{alpha
  vertex bound}) then gives a quasiconvexity bound for $\{[a_i]\}$, 
completing the proof of Theorem \ref{Curve replacements}.
\end{proof}

\subsection*{Proof of the main theorem}
The proof of Theorem \ref{theorem:handlebody} can now be easily
assembled: 

Given two vertices of $\Delta(M,S)$, represent them in minimal
position as $a,b$ and let $J_0\subset a$ be a good interval containing
$a\intersect b$  (by Lemma \ref{good intervals}).
Proposition \ref{prop:nested} gives us a nested curve replacement
starting with $(a_0,J_0)$, which is well-nested by choice of $J_0$,
and whose vertices are all in $\Delta(M,S)$.
By Theorem \ref{Curve replacements}, this set
is $K$-quasiconvex. It follows that $\Delta(M,S)$ is $K$-quasiconvex.
\qed


\begin{thebibliography}{10}

\bibitem{bowditch:complex}
B.~Bowditch, \emph{Intersection numbers and the hyperbolicity of the curve
  complex}, Preprint, Southampton.

\bibitem{brock-canary-minsky:ELCII}
J.~Brock, R.~Canary, and Y.~Minsky, \emph{Classification of {Kleinian} surface
  groups {II}: the ending lamination conjecture}, in preparation.

\bibitem{casson-gordon:reducing}
A.~J. Casson and C.~McA. Gordon, \emph{Reducing {H}eegaard splittings},
  Topology Appl. \textbf{27} (1987), no.~3, 275--283.

\bibitem{hempel:complex}
J.~Hempel, \emph{3-manifolds as viewed from the curve complex}, Topology
  \textbf{40} (2001), no.~3, 631--657.

\bibitem{johannson:3mfds}
K.~Johannson, \emph{Topology and combinatorics of 3-manifolds}, Lecture Notes
  in Mathematics, vol. 1599, Springer-Verlag, Berlin, 1995.

\bibitem{masur-minsky:complex1}
H.~A. Masur and Y.~Minsky, \emph{Geometry of the complex of curves {I}:
  Hyperbolicity}, Invent. Math. \textbf{138} (1999), 103--149.

\bibitem{masur-minsky:complex2}
\bysame, \emph{Geometry of the complex of curves {II}: Hierarchical structure},
  Geom. Funct. Anal. \textbf{10} (2000), 902--974.

\bibitem{minsky:boundgeom}
Y.~Minsky, \emph{Bounded geometry in {Kleinian} groups}, Invent. Math.
  \textbf{146} (2001), 143--192.

\bibitem{minsky:ELCI}
\bysame, \emph{Classification of {Kleinian} surface groups {I}: models and
  bounds}, preprint, 2002.

\bibitem{penner-harer}
R.~Penner and J.~Harer, \emph{Combinatorics of train tracks}, Annals of Math.
  Studies no. 125, Princeton University Press, 1992.

\bibitem{schleimer:distance2}
S.~Schleimer, \emph{The disjoint curve property}, preprint, 2002.

\end{thebibliography}

\providecommand{\bysame}{\leavevmode\hbox to3em{\hrulefill}\thinspace}

\end{document}